\documentstyle[rkchngs,12pt]{article}
\def\half{{1 \over 2}}
\def\bigspace{\bigskip \\}
\def\medspace{\medskip \\}
\def\smallspace{\smallskip \\}

\def\supsetn{\raisebox{-.8ex}[0in]{$\:\prime$}\kern-.50em$\supseteq$}

\def\rar{\rightarrow}

\font\gothic=eufm10 at 12pt%
\def\g{{\gothic g}}

\font\bbten=msbm10
\font\bbseven=msbm7
\newfam\bbfam \textfont\bbfam=\bbten \scriptfont\bbfam=\bbseven
\def\bb{\fam\bbfam\bbten}
\newcommand{\bbr}{{\bb R}}

\def\pd#1{{\partial \over \partial {#1}}}
\def\e6{$\{e_1,\ldots,e_6\}$}
\def\w6{$\{\omega_1,\ldots,\omega_6\}$}
\def\x6{$(x_1(t),\ldots,x_6(t))$}

\def\gm{\gamma}
\def\lm{\lambda}
\def\al{\alpha}
\def\om{\omega}
\def\tt{\theta}

\def\er{$\{e_1,\ldots,e_r\}$}
\def\en{$\{e_1,\ldots,e_n\}$}
\def\tetan{$\{\theta_1,\ldots,\theta_n\}$}

\newtheorem{theorem}{Theorem:}[section]
\newtheorem{lemma}[theorem]{Lemma:}
\newtheorem{prop}[theorem]{Proposition:}

\newtheorem{df}[theorem]{Definition:}
\newtheorem{formula}[theorem]{}

\def\remarks{{\bf Remarks \ }}
\def\remark{{\bf Remark \ }}
\def\example{{\bf Example: \ }}
\def\Bproof{{\bf Proof:\ }}
\def\Eproof{\ $\Box$}
\newenvironment{proof}{{\bf Proof:\ }}{\ $\Box$}

\def\IMSmarkvadjust{0 pt}
\def\IMSmarkhadjust{0 pt}
\def\IMSmarkhpadding{0 pt}
\def\IMSpubltext{Published in modified form:}
\def\SBIMSMark#1#2#3{
 \font\SBF=cmss10 at 10 true pt
 \font\SBI=cmssi10 at 10 true pt
 \setbox0=\hbox{\SBF \hbox to \IMSmarkhpadding{\relax}
                Stony Brook IMS Preprint \##1}
 \setbox2=\hbox to \wd0{\hfil \SBI #2}
 \setbox4=\hbox to \wd0{\hfil \SBI #3}
 \setbox6=\hbox to \wd0{\hss
             \vbox{\hsize=\wd0 \parskip=0pt \baselineskip=10 true pt
                   \copy0 \break%
                   \copy2 \break%
                   \copy4 \break}}
 \dimen0=\ht6   \advance\dimen0 by \vsize \advance\dimen0 by 8 true pt
                \advance\dimen0 by -\pagetotal
	        \advance\dimen0 by \IMSmarkvadjust
 \dimen2=\hsize \advance\dimen2 by .25 true in
	        \advance\dimen2 by \IMSmarkhadjust

  \openin2=publishd.tex
  \ifeof2\setbox0=\hbox to 0pt{}
  \else 
     \setbox0=\hbox to 3.1 true in{
                \vbox to \ht6{\hsize=3 true in \parskip=0pt  \noindent  
                {\SBI \IMSpubltext}\hfil\break
                {\it J. Dynam. Control Sys.}~{\bf 1} (1995), 535--549
 
                \vfill}}
  \fi
  \closein2
  \ht0=0pt \dp0=0pt
 \ht6=0pt \dp6=0pt
 \setbox8=\vbox to \dimen0{\vfill \hbox to \dimen2{\copy0 \hss \copy6}}
 \ht8=0pt \dp8=0pt \wd8=0pt
 \copy8
 \message{*** Stony Brook IMS Preprint #1, #2. #3 ***}
}

\begin{document}
\title{A Note on Carnot Geodesics in Nilpotent \\ Lie Groups}
 
\author{Christophe Gol\'e and Ron Karidi
\thanks{
The first author was partially supported by an NSF Postdoctoral
Fellowship, and the second author was
partially supported by a Rothschild Postdoctoral Fellowship.}}
 
\date{}
 
\maketitle
%

\SBIMSMark{1995/6}{May 1995}{}

\begin{abstract}
We show that strictly abnormal geodesics arise in graded nilpotent
Lie groups. We construct such a group, 
for which some Carnot geodesics are strictly abnormal and, 
in fact,  not normal in any subgroup.
In the 2-step case we also prove that these geodesics are always smooth.
Our main technique is based on 
the equations for the normal and abnormal 
curves, that we derive (for any Lie group)
explicitly in terms of the structure constants. 
\end{abstract}

\thispagestyle{empty}

\section{Introduction} 

Our work is motivated by the problem of differentiability of Carnot
geodesics (or minimizers) in nilpotent Lie groups.

A Carnot (or sub-Riemannian) structure on a manifold $G$ is given 
by a smoothly varying 
distribution $D$ (i.e. a field of tangent subspaces) and a smoothly 
varying inner product on  this distribution.
A {\it horizontal curve} is an absolutely continuous curve on $G$ 
which is tangent to $D$ wherever it is differentiable.
The inner product on $D$ enables one to define the length of 
a horizontal curve, and it is then natural to study Carnot geodesics, 
or minimizers for this length. 
A {\it minimizer} is an absolutely continuous horizontal curve $x$ in 
$G$ which is such that,  for each $t$, 
there exists $\epsilon >0$ such that $x$ minimizes the length between
$x(t_0)$ and $x(t_1)$ whenever $t_0, t_1$ are in $(t-\epsilon, t+\epsilon)$.

Until recently, 
it was not understood that minimizers could be of two different, {\it but non mutually exclusive} types.
%
One type is given by the projections of solutions of a 
Hamiltonian system (see Section \ref{Snormal}), 
which is in a sense the
Legendre transform of the inner product on $D$. 
This generalizes the Riemannian situation. 
Such curves are called {\it normal}. Normal
curves are known to be differentiable minimizers \cite{Hamen}.
\\
The other type of minimizer belongs to a category of horizontal curves 
called {\it abnormal} 
or {\it singular}. Although originally given by the Maximum Principle, 
Hsu \cite{Hsu} (see also \cite{Mont1})
shows that they are projections onto $G$ of characteristic curves 
(in the symplectic sense) of the annihilator of $D$ in $T^*G$. 
See Definition \ref{Dabnormal}.


Contrary to the normal curves, the abnormal ones need not be minimizing.
If they are minimizing, they are called
{\it abnormal minimizers}.
That such curves exist was only proven recently by Montgomery
\cite{Mont2}.
For rank-2 distributions,
Liu and Sussmann \cite{LS}
point to a generic class of abnormal curves that also minimize.
They are called {\it regular abnormal extremals}, and are used in 
Section \ref{Sreg}.
These curves were also studied by Bryant and Hsu in their
paper on  rigid curves \cite{BH}.
Again, let us emphasize that 
the same curve can be  both normal and abnormal: 
it can have a lift in the cotangent bundle that makes it abnormal and
another one that makes it normal. If a curve is abnormal but not normal,
we say that it is {\it strictly abnormal.}

In this paper we study Carnot minimizers in nilpotent Lie groups. 
This context is important 
because, under natural conditions, any Carnot manifold, 
viewed as a metric space,
asymptotically looks like some nilpotent Lie group.  
More precisely,
given any metric space $(M,d)$, Gromov (\cite{Gromov} Section 1.4.B) defines
the tangent cone to $M$ at $q \in M$ to be the Hausdorff
limit of the metric spaces (M, $\lambda d$), with base point at $q$,
when $\lambda \rar \infty$.
Mitchell \cite{Mitchell} proves that if a Carnot manifold $(M,D)$, 
is regular
at $q$, then the tangent cone at $q$ is isometric (as a metric space) to 
a Carnot graded nilpotent Lie group (called the nilpotentization), 
defined using the bracket relations on $T M$ 
(see Section \ref{Seqns} below).

Our most surprising result (see Theorem \ref{Tcountereg})
is the construction 
of a graded nilpotent Lie group, where we find
many abnormal minimizers that are not normal, i.e. we find
{\it strictly abnormal minimizers}.
In light of Mitchell's theorem, this shows that 
the analysis of abnormal minimizers cannot be avoided
in general Carnot spaces.

That strictly abnormal minimizers exist is relevant to 
a natural problem posed by Hamen- st\"adt:
{\it are Carnot minimizers always differentiable ?}
(see \cite{Hamen},\cite{Mont2}). Indeed, since normal minimizers are
automatically differentiable, Hamenst\"adt's question reduces to: {\it are abnormal
minimizers always differentiable?}
There are examples of non differentiable
abnormal curves, even in the nilpotent Lie group situation studied here, 
but none that are minimizers.
To solve this problem in Lie groups, Hamenst\"adt suggested
to try to prove that any minimizer is normal in a subgroup,
say $H < G$, with the Carnot structure given by the distribution
$D_H=D\cap TH$ and the restricted inner product. 
The existence of such an $H$ for each minimizer would obviously imply that 
any minimizer is differentiable.
This program was carried out successfully by Montgomery in \cite{Mont1}
in the case where
$G$ is a compact connected Lie group and $D$ is the left invariant 
distribution orthogonal to its maximal torus. 

In contrast, the strictly abnormal minimizers that we exhibit in
Theorem \ref{Tcountereg}, 
are not normal in any subgroup.
However, we can carry out Hamenst\"adt's 
program in the case where $G$ is a 2--step
nilpotent Lie group (i.e.  its Lie algebra \g \ satisfies 
[\g , [\g ,\g ]]$=0$) and $D$ is a left invariant distribution such that 
$D\oplus$[\g ,\g ]=\g \ (Theorem \ref{Tsmooth}).
Note that the latter condition is satisfied for graded nilpotent Lie algebras.
This extends known results on the so-called Gaveau-Brockett problem
\cite{Br},\cite{Gav}.

Our methods rely on deriving the equations for the normal and abnormal curves 
purely in terms of the structure constants.
These equations appear in Section \ref{Seqns}. 
We also found the methods in \cite{Ka} useful in our investigations.

We are very grateful to Richard Montgomery for many useful discussions.
The first author would like to thank Stanford University for its hospitality
during the time when part of this work was done.

\section{Abnormal and Normal equations}
\label{Seqns}

Let $G$ be a Lie group, and $D$ a left invariant distribution.
We identify D with  a left translation of a subspace of the Lie algebra \g, 
that we will also denote by $D$.
Choose a left invariant frame \er\  for $D$,
and complete it to a basis \en \ of the whole Lie algebra.
We give $D$ a metric that makes \er  an orthonormal basis.
Let \tetan \ be the dual co-frame to \en.
We write a vector field in $TG$ as $\sum _1^n\gm_i e_i$, \
where $\gm_1,\ldots,\gm_n$ are the coordinate functions on the fiber of $TG$.
Likewise, a covector is written $\sum_1^n \lm_i \theta_i$.
A vector  in $D$ (resp. a covector  in $D^\perp$)
can be written as 
$\sum_1^r \gm_i e_i$ \ (resp. $\sum_{r+1}^n \lm_i \theta_i$).
\medspace
The structure constants $\al_{ijk}$ of \g \ with respect to the basis 
\en \ are defined by:
$$[e_i,e_j]=\sum_{k=1}^n \al_{ijk} e_k \ . $$

Note that $\alpha_{ijk}=-\alpha_{jik}$. 
Since $e_i$ is the Hamiltonian vector field for $\lm_i$, we also have:
$$\{\lm_i,\lm_j\}=-\sum_{k=1}^n \al_{ijk} \lm_k \ , $$
with respect to the standard Poisson brackets of functions on \g$^*$. 
%

Denote $D^1=D, \
D^{i+1}=D^i + [D,D^i]$
\ (if this sum is direct, one says that \g \ is {\it graded}).
If there exists $r$ for which $D^r=$ \g \  we say that $D$ 
is {\it bracket-generating.}
Define also $V^1=D^1 \ , \ V^i=D^i / D^{i-1}$ and
$Gr $ \g $=V^1 \oplus \ldots \oplus V^r$. The latter is a graded nilpotent
Lie algebra and the associated simply connected Lie group is called
the {\it nilpotentization} of $G$, which is also endowed
with a Carnot metric.
\\
\remark
In the case of a general Carnot manifold $G$, the above makes
sense locally, at a point $q \in G$, if one assumes, in addition, 
that $r(q)$ is locally constant at that point. 
The theorem of Mitchell \cite{Mitchell} alluded
to in the introduction relates the local metric properties of a Carnot manifold
with those of its nilpotentization.
%
%
\subsection{The Abnormal Equations}
\label{Sabnormal}
We begin by giving a rigorous definition of abnormal curves and minimizers:
\begin{df}
\label{Dabnormal} An abnormal curve is a horizontal curve which is the 
projection onto $G$ of an absolutely continuous curve in the annihilator 
$D^\perp \subset T^*G$ of
$D$, with square integrable derivative, which does not intersect the zero 
section and whose derivative, whenever it exists, is in  
the kernel of the canonical symplectic form restricted to  $D^\perp$. 
An abnormal minimizer is an abnormal
curve which is a minimizer, in the sense given in the introduction. A strictly
abnormal curve (resp. minimizer) is an abnormal curve (resp. minimizer) which 
is not normal, in the sense of Definition \ref{Dnormal}.
\end{df}

This definition makes it clear that being abnormal is independent
of the parameterization. In \cite{Mont1}, Proposition 1, Montgomery shows 
that the above definition of abnormal curve is equivalent to three other ones, 
which we will not use in this paper,
but will state for the further confusion of the reader. 
Let $x$ be a curve in $G$ and  $\zeta$ be 
an absolutely continuous curve in $T^*G$ which does not intersect the 
zero section and whose derivative is square integrable. 
Suppose that $x=\pi(\zeta)$ and $x$ and $ \zeta$ satisfy the above definition. 
Then, the following are equivalent to the above definition.
\begin{enumerate}
\item  $x$ is an abnormal extremal in the sense of the
Pontryagin maximum principle of control theory (this is the original definition
of abnormal curves).

\item  $\zeta$ annihilates the image
of the differential $d (end(x(t)))$ at each $t$, where $end$ is the map associating
to a curve its endpoint.

\item $x$ is horizontal, 
$\zeta \in D^\perp$ and $\zeta(t)= (D\Phi_t^T)^{-1}\zeta(0)$ 
where $\Phi_t$ is any time dependent flow which generates the curve $x$.
\end{enumerate}

We now follow the derivation of the abnormal equations in (\cite{Mont1}
Section 4).
Remember that the canonical 1-form on $T^*G$, call it $\eta$, is defined 
by $\eta_\alpha (v)=\alpha( \pi_* v)$, where $\alpha \in T^*G$ is the base point
of the vector $v \in T_\alpha (T^*G)$.
We claim that $\eta_\alpha=\sum_{i=1}^n \lm_i(\alpha) \tt_i$, 
where $\tt_i$ are viewed as 1-forms on $T^*G$.
Let $v=\sum \gm_i e_i + \sum h_i \pd{\lm_i} \in T(T^*G)$, then
by the definition of $\eta$:
$$ 
\eta_\alpha(v)=\alpha (\pi_* v)=
\alpha(\sum \gm_i e_i)=
\sum \lm_i(\alpha) \tt_i (\sum \gm_i e_i)=\sum \lm_i(\alpha) \gm_i  \ , $$
which coincides with $\sum \lm_i(\alpha) \tt_i(v)$.
In particular, $\eta_{|_{D^\perp}}=\sum_{i=r+1}^n \lm_i(\alpha) \tt_i$, and,
since $\om=d \eta$,
$$\omega_{|_{D^\perp}}= 
\sum_{i=r+1}^n \ (d \lm_i \wedge \tt_i + \lm_i d \tt_i) \ . $$
\vspace{.2cm}

Let $(x(t),\lm(t)) \in T^*G$ be such that $x(t)$ is an abnormal curve.
Then
$$ x'={dx\over dt}= \sum_{i=1}^r \gm_i e_i \ \ , \ \ 
\lm'={d\lm\over dt}=
\sum_{i=r+1}^n {\lm_i}' \pd{\lm_i} \ \ , \ \ \lm \ne 0 \ \ , $$ 
and
\begin{formula}
\label{Fomega}
$$
0=\omega_{(x,\lm)} \ ((x',\lm'),\cdot)=
\sum_{i=r+1}^n ({\lm_i}'  \tt_i - \gm_i d \lm_i) + 
\sum_{k=r+1}^n \ \lm_k d \tt_k((x',\lm'),\cdot) \ \ .
$$
\end{formula}
The Maurer-Cartan equations are:
$$
d \theta_k + \half \sum_{i,j=1}^n \al_{ijk} \theta_i \wedge \theta_j = 0 \ . 
$$
Therefore 
$$d \tt_k((x',\lm'),\cdot)=
-\half \sum_{i,j=1}^n \alpha_{ijk} (\gm_i \tt_j-\gm_j \tt_i)= 
\sum_{i,j=1}^n \alpha_{ijk} \gm_j \tt_i \ . $$
Plug this in (\ref{Fomega}), using $\gm_i=0$ for $i>r$:
$$ 0=\sum_{i=r+1}^n \ ( \ {\lm_i}' + \sum_{j=1}^n \sum_{k=r+1}^n 
\alpha_{ijk} \gm_j \ ) \ \tt_i  
+ \sum_{i=1}^r \ ( \ \sum_{j=1}^n \sum_{k=r+1}^n \alpha_{ijk} \gm_j \ ) \ \tt_i
\ . $$
Since the $\tt_i$'s are linearly independent we get 

\begin{formula}
\label{Eabnormal}
The abnormal equations: 
$$
\sum_{j=1}^r \sum_{k=r+1}^n \al_{ijk} \gm_j \lm_k = 0 \ ,
\ \ \ {\rm for} \ \ i=1,\ldots,r \ .
$$
\
$${\lm_i}' + \sum_{j=1}^r \sum_{k=r+1}^n \al_{ijk} \gm_j \lm_k = 0 \ , 
\ \ \ {\rm for} \ \ i=r+1,\ldots, n \ . 
$$
\
$$
\gm_{r+1}= \ldots=\gm_n=0 
$$
$$
\lm_1=\ldots=\lm_r=0 $$
\end{formula}
\remark
These equations are mixed algebraic-differential equations.
Unlike the normal equations 
(which, as we will see, are defined by a single vector field),
the abnormal ones cannot be expressed as ordinary differential equations.
Note also that the notion of abnormal curve only depends on $D$, 
and not on the metric.
%
\subsection{The Normal Equations}
\label{Snormal}

To distinguish the cotangent lifts of normal and abnormal curves,
we denote the covector frame coordinates in the normal case by 
$h_1,\ldots,h_n$ (instead of $\lm_1,\ldots,\lm_n$).

\begin{df}
\label{Dnormal}
A normal curve is the projection onto $G$ of a
solution of the Hamiltonian system in $T^*G$ with Hamiltonian 
$H(x,h)=\half \sum_1^r h_i^2$. 
\end{df}
It is known \cite{Stric} that 
normal curves are smooth minimizers.
We can write the Hamiltonian
vector field in the frame coordinates as $(\gm, h')$. 
It will satisfy:

\begin{formula}
\label{Enormal}
The normal equations:
$$
\begin{array}{l}
\gm_i=h_i \ \ \ {\rm for} \ \ i=1,\ldots,r \\
\gm_j=0 \ \ \ \ {\rm for} \ \ j=r+1,\ldots,n
\end{array}
$$

$${h_i}'+\sum_{j=1}^r \sum_{k=1}^n \al_{ijk} h_j h_k = 0 
\ \ \ {\rm for} \ \ i=1,\ldots,n \ . $$
\end{formula}

\ \bigspace
The first two equations are given by the fact that, in a canonical
system of coordinates $(q,p)$ on $T^*G$, $h_i$ can be
seen as the fiber linear function:
$$
h_i(q,p)=p(e_i(q)),\  q\in G,\  p\in T^*_qG
$$
To get the third equation observe that for any function $f$ on $T^*G$
and a solution $z(t)$ to the Hamiltonian equations:
$$ {d \over d t} (f(z(t)) = \{f,H\} \ . $$
In particular 
${h_i}'=\{h_i,H\} \ \ \ {\rm for} \ \ i=1,\ldots,n$ \ and thus 

$${h_i}'=\{h_i,\half \sum_1^r h_j^2\}=
\half \sum_1^r \{h_i,h_j\} h_j = 
- \sum_{j=1}^r \sum_{k=1}^n \al_{ijk} h_j h_k.$$
\bigspace

Consider now a {\it nilpotent} Lie group $G$ and 
assume that the distribution D is complementary to [\g,\g],
i.e. $D \oplus [$\g,\g$]=$\g.
Note that this condition is satisfied when $G$ is graded
and $D$ is bracket-generating.
In fact, the meaning of this assumption is
that $D$ is minimal in the sense that
no proper subspace
of $D$ will bracket generate the full Lie algebra.
This also means that $\alpha_{ijk}=0 $ for $k=1,\ldots,r$ and thus the second normal 
equation reduces to:

$${h_i}'+\sum_{j=1}^r \sum_{k=r+1}^n \al_{ijk} h_j h_k = 0 \ . $$

\begin{prop}
Under these assumptions, the one-parameter horizontal subgroups
are normal curves for any left invariant metric on D.
\end{prop}
\Bproof
Let $x(t)$ be a ``left invariant'' horizontal curve, i.e.
$x'(t)=\sum \gm_i e_i$, where $\gm_i$ are constants and
$\gm_{r+1}=\ldots=\gm_n=0$.
We choose the covector part to have the following coordinates:
\\
\begin{eqnarray*}
& h_i=\gm_i & i=1,\ldots,r 
\\
& h_i=0 & i=r+1,\ldots,n
\end{eqnarray*}

It is then easy to verify that the normal equations hold.
\Eproof
\medspace
\remarks
\begin{enumerate}
\item
Note that one can always assume that a normal curve $x$ is parameterized by
arc-length. Indeed, $\half {\Vert {x'} \Vert}^2 =H$ is constant.
\item
As we have said in the introduction, being abnormal and being normal are not mutually exclusive properties.
The same curve $x(t)$  in $G$ can have lifts $(x,\lm)$ and $(x,h)$,
with the first making $x$ abnormal, the other normal.
When that happens, one can choose the arc-length parameterization 
for both.
\end{enumerate}
\example
Take the Engel algebra, $[e_1,e_2]=e_3 \ , \ [e_1,e_3]=e_4$,
with $D=Span\{e_1,e_2\}$. Then the only abnormal curves are
tangent to $e_2$, i.e. they are one-parametric subgroups, and by the above
proposition, they are also normal.

\subsection{Regular abnormal extremals}
\label{Sreg}

Liu and Sussmann \cite{LS} show that if $D$ is a two dimensional
distribution in a Carnot manifold M of dimension $\ge 3$, there is
an efficient way of finding lots of abnormal geodesics, or minimizers.

Namely, they introduce the notion of {\it regular abnormal extremals}
and prove that all such regular abnormal extremals are in fact minimizers
(\cite{LS}, Theorem 5).  
The following definition is equivalent to theirs, and can be extracted
from their Proposition 6 and the beginning of their Section 6.2.
As before, let $D^k$ be the set of Lie brackets of order $k$ or 
less of vector fields in $D$,
and $(D^k)^\perp$ the annihilator of this set in the cotangent bundle.

\begin{df}
\label{regular}
A curve $x(t)$ in $G$, parameterized by arclength, is called a regular 
abnormal extremal (or regular abnormal, in short) if it has a lift
$(x(t), \lambda(t))$ which satisfies the abnormal equations \ref{Eabnormal}
and such that
$$ \lambda(t)\in (D^2)^\perp-(D^3)^\perp.$$
\end{df}

Another interesting property of regular abnormals, which 
we will not use here,
is that they are projections of  integral curves of a certain vector
field $\chi(D)$ in $ (D^2)^\perp-(D^3)^\perp$. This is in fact the 
definition of a regular abnormal extremal given in \cite{LS}. 
Liu and Sussman have also genericity results
showing that, roughly, among all lifts of abnormal curves 
parameterized by arclength, the regular
abnormal extremals are prevalent.

\section{Strictly abnormal geodesics in nilpotent Lie groups}

In this section, we produce examples of strictly abnormal curves 
in a graded nilpotent Lie group $G$. 
Note that such examples exist in non-nilpotent Lie groups
(see for example \cite{LS}, Section 9.5).
Our examples also have the property that
they are not normal in any proper subgroup of $G$.

Let \g \ be the 6 dimensional real Lie algebra
spanned by  \e6 with the following relations:
$$ [e_1,e_2]=e_3 \ , \ [e_1,e_3]=e_4 \ , \ 
   [e_2,e_3]=e_5 \ , \ [e_1,e_4]=e_6 \ .
$$
Let the distribution $D=Span\{e_1,e_2\}$. 
It satisfies $D\oplus [$\g$,$\g$]=$\g.
We put on $G$
the Carnot metric that makes the frame $e_1,e_2$ of the distribution $D$
orthonormal. We are looking for an abnormal curve, say \x6, 
which cannot be normal.
Writing $x'=\gm_1 e_1 + \gm_2 e_2$,
we first find a lift of $x$ to the cotangent bundle,
which satisfies the abnormal equations.
Denote (as in Section \ref{Seqns}) 
the cotangent part by $\sum \lm_i(t) \theta_i$ where 
$(\theta_1,\ldots \theta_6)$ 
is a basis of left invariant 1-forms on $G$ dual to \e6.

The abnormal equations given in (\ref{Eabnormal}) become in this case:
\begin{eqnarray*}
\gm_2 \lm_3=0
\\
-\gm_1 \lm_3=0
\\
\lm_3'-\gm_1 \lm_4 - \gm_2 \lm_5=0 
\\
{\lm_4}' - \gm_1 \lm_6=0
\\
{\lm_5}'=0
\\
{\lm_6}'=0
\end{eqnarray*}
Remember that, for an abnormal, $\lm_1=\lm_2=0$ and note that here $\lm_3=0$ as
well (otherwise $x$ would be the trivial, constant solution), 
so that the $\lm_3'$ term in the third equation is actually zero.
We will  assume that $x$ is parameterized by arc-length, i.e. 
$\gm_1^2+\gm_2^2=1$.
For our class of examples, 
we will seek a solution $x$ with $\lm_5,\lm_6 \ne 0$.
\bigspace
From $\gm_1={{\lm_4}' \over \lm_6}  \ , \ 
\gm_2=- \gm_1 {\lm_4 \over \lm_5} = 
- {\lm_4 {\lm_4}' \over \lm_5 \lm_6}$, and 
$\gm_1^2+\gm_2^2=1$ we get an o.d.e for $\lm_4$:

\begin{formula}
\label{Emain}
$$ ({\lm_4}')^2 (\lm_4^2+\lm_5^2) = \lm_5^2 \ \lm_6^2 \ . $$
\end{formula}

\begin{theorem}
\label{Tcountereg}
Let $\lm_4$ be a solution of equation (\ref{Emain}),
with $\lm_5,\lm_6 \ne 0$. Then any solution to the
time dependent o.d.e 
$x'=\gm_1e_1+\gm_2e_2$ with $x(0)=0$, where $\gm_1, \gm_2$ satisfies
$\gm_1={{\lm_4}' / \lm_6}
\ , \ \gm_2=-{\lm_4 {\lm_4}' / \lm_5 \lm_6}$ 
is a strictly abnormal minimizer. Such a curve cannot be normal
in any subgroup of $G$ either.
\end{theorem}
\begin{proof}
That such an $x$ is abnormal derives directly from the abnormal equations
and our discussion above. To see that it is a minimizer we 
observe that it is a regular abnormal extremal.
Indeed, $\lm$ belongs to $(D^2)^\perp$ since
$(x,\lm)$ satisfies the abnormal equations, and
$\lm_5 \ne 0$
implies that $\lm \not \in (D^3)^\perp$.

We now argue that 
$\gm_1,\gm_2$ cannot be constant. By (\ref{Emain}),
$\lm_4'\neq 0$. But then, the same equation tells us that, since $\lm_4$
is not constant, neither is $\lm_4'$.
From that it easily follows that $\gm_1$ and $\gm_2$ are not constant.
The following proposition tells us
that, because of this, $x$ cannot be normal in $G$, 
hence it is a strictly abnormal minimizer.

Finally, the fact that $\gm_1$ and $\gm_2$ are non constant also tells us that
$x$ cannot be embedded in any proper subgroup of $G$. If it were,
the pull-backs of the tangent 
vectors to $x$ would all belong to a proper subalgebra.
But since $\gm_1,\gm_2$ are not constant, there exist two pull-backs 
that span $D$ as a vector space,
and hence generate \g \ as a Lie algebra.
\end{proof}
 
\begin{prop}
\label{Pleftinv}
The only normal abnormals 
in this Lie group are the left invariant
curves, i.e. integral curves of the left invariant vector fields
(for which $\gm_1$ and $\gm_2$ are constant).
\end{prop}
\begin{proof}
Let $x$ be an abnormal curve which is also normal.
As before, we assume that $x$ is parameterized by arc-length.
This implies the existence of two cotangent lifts: \  
$(x,\lm)$ and $(x,h)$ for which the corresponding 
$(\gm,\lm)$ and $(\gm,h)$ satisfy the abnormal and normal equations
respectively.

We first prove that if $\lm_5=0$ or $\lm_6=0$,  
then $x$ is left invariant.
Looking at (\ref{Emain}), we see that if one of $\lm_5$ or $\lm_6$
is zero, $\lm_4'$ must also be zero, so $\lm_4$ is constant.
This gives the linear equation: $\gm_1 \lm_4 + \gm_2 \lm_5 = 0$,
with constant coefficients. 
Since $\lm_4 \ne 0$, the assumption of arc-length
parameterization implies that $\gm_1$ and $\gm_2$ must be constant.
%
We now assume that $\lm_5,\lm_6 \ne 0$ and get a contradiction.

We proceed by combining the abnormal equations and 
\begin{formula}
\label{Exnormal}
The normal equations:
\bigspace
$ (a) \ \ \ {\gm_1}'+\gm_2 h_3=0 \\
(b) \ \ \ {\gm_2}'-\gm_1 h_3=0 \\
(c) \ \ \ {h_3}'-\gm_1 h_4 - \gm_2 h_5=0 \\
(d) \ \ \ {h_4}'-\gm_1 h_6=0 \\
(e) \ \ \ {h_5}'=0 \\
(f) \ \ \ {h_6}'=0 
$
\end{formula}
We start by deriving $h_3$.
Differentiating (\ref{Emain}) we get 
\begin{formula}
\label{Flm4pp}
$${\lm_4}'' = - {\lm_4 ({\lm_4}')^4 / ({\lm_5}^2 {\lm_6}^2}) \ , $$
\end{formula}
which gives, after differentiating $\gm_2=- \lm_4 {\lm_4}' / \lm_5 \lm_6$ 
$$ {\gm_2}'=-{({\gm_4}')^2 \over \lm_5^3 \lm_6^3} 
(\lm_4^2 ({\lm_4}')^2 - \lm_5^2 \lm_6^2) \ \ \ \ 
{\rm and \ by \ (\ref{Emain})} \ \ 
{\gm_2}'=-{({\lm_4}')^4 \over \lm_5 \lm_6^3} \ . $$
Now, $h_3={\gm_2}' / \gm_1$, and $\gm_1={\lm_4}'/\lm_6$. Hence
$$h_3=-{({\lm_4}')^3 \over (\lm_5 \lm_6^2)} $$ 
(recall that we can divide by ${\lm_4}'$ which is never zero
by \ref{Emain}).
\medspace
Differentiate this equality and substitute in (\ref{Exnormal}c) to get
$$-{3 ({\lm_4}')^2 {\lm_4}'' \over \lm_5 {\lm_6}^2} =
{{\lm_4}' \over \lm_6} h_4 - {\lm_4 {\lm_4}' \over \lm_5 \lm_6} h_5 \ .$$
Together with (\ref{Flm4pp}), it gives
$$ 3 \lm_4 ({\lm_4}')^5 / {\lm_5}^2 {\lm_6}^3 =
\lm_5 h_4 - \lm_4 h_5 \ . $$
\medspace
Now, ${h_4}'=\gm_1 h_6=
(h_6 / \lm_6){\lm_4}'$ so $h_4=(h_6 / \lm_6) \lm_4 + c_4$,
where $c_4$ is a constant.
This implies 
$$ {3 \lm_4 ({\lm_4}')^5  \over {\lm_5}^2 {\lm_6}^2} =
c_4 \lm_5 \lm_6 + \lm_4 (\lm_5 h_6 - \lm_6 h_5) $$

We introduce generic constants $\alpha_1, \alpha_2, \alpha_3$, which 
we will use liberally to denote any constants (the actual
value they represent may change from one equation to the other).
The above equation writes as:
$$ \lm_4 (\al_1 + \al_2 ({\lm_4}')^5)=\al_3 
\ \ \ \ \ (\al_2 \ne 0) \ . $$
We differentiate, divide by ${\lm_4}'$ and use (\ref{Flm4pp})
to get the following equation:
$$\al_1 + \al_2 ({\lm_4}')^5 + \al_3 \lm_4^2 ({\lm_4}')^7 = 0 
\ \ \ \ (\al_2,\al_3 \ne 0) \ . $$
Or \ $\al_1 + ({\lm_4}')^5 (\al_2+\al_3 {\lm_4}^2 ({\lm_4}')^2) = 0 
\ \ \ \ \ (\al_2,\al_3 \ne 0) \ . $
\bigspace
From (\ref{Emain}) we see that 
${\lm_4}^2 ({\lm_4}')^2 = {\lm_5}^2 ({\lm_6}^2 - ({\lm_4}')^2)$,
which gives
$$\al_1 + \al_2 ({\lm_4}')^5 + \al_3 ({\lm_4}')^7 = 0 
\ \ \ \ \ (\al_3 \ne 0) \ . $$
Hence ${\lm_4}'$ must take a finite set of values, but by (\ref{Emain})
so does $\lm_4$, which is continuous. So $\lm_4$ is constant, and
${\lm_4}'=0$. Contradiction.
\ \end{proof}
\bigspace

This phenomena, where the only normal abnormals are left invariant
curves, is not rare. In fact it occurs in many examples.
However, one can build counter examples
to this, in rank-2 distributions,
whenever there exists a normal curve with a constant $h_3$.
It is not hard to see in this case that, taking $\lm_3=0$ and
$\lm_i=h_i$ for $i>3$, gives a solution to the abnormal equations.

As an example we take the free nilpotent Lie algebra of step 4
on two generators.
This Lie algebra is of dimension 8, with the following relations:
$$ [e_1,e_2]=e_3 \ , \ [e_1,e_3]=e_4 \ , \ 
   [e_2,e_3]=e_5 \ , \ [e_1,e_4]=e_6 \ , \
$$
$$
   [e_1,e_5]=e_7 \ , \ [e_2,e_4]=e_7 \ , \
   [e_2,e_5]=e_8 \ .
$$
Let the distribution $D=Span\{e_1,e_2\}$.
The abnormal equations become:
\begin{eqnarray*}
\gm_2 \lm_3=0
\\
-\gm_1 \lm_3=0
\\
-\gm_1 \lm_4 - \gm_2 \lm_5=0 
\\
{\lm_4}' - \gm_1 \lm_6 - \gm_2 \lm_7=0
\\
{\lm_5}' - \gm_1 \lm_7 - \gm_2 \lm_8=0
\\
{\lm_6}'=0
\\
{\lm_7}'=0
\\
{\lm_8}'=0
\end{eqnarray*}
And the normal equations are:
\begin{eqnarray*}
{\gm_1}'+\gm_2 h_3=0 \\
{\gm_2}'-\gm_1 h_3=0  \\
{h_3}'=\gm_1 h_4 + \gm_2 h_5 \\
{h_4}'=\gm_1 h_6 + \gm_2 h_7 \\
{h_5}'=\gm_1 h_7 + \gm_2 h_8 \\ 
{h_6}'=0 \\ 
{h_7}'=0 \\ 
{h_8}'=0 
\end{eqnarray*}
\\
Any curve given by integrating
$$ x'(t)= (-\sin t) e_1 + (\cos t) e_2,$$
is, on one hand, not left invariant ($\gm_1$ and $\gm_2$
are not constant), and on the other hand  both normal 
and abnormal. 
As an abnormal lift of $x$ to the cotangent bundle we take: 
$$ \lm_1=\lm_2=\lm_3=0 \ , \\
\lm_4=\cos t \ , \ \lm_5=\sin t \ , \ 
\lm_6=1 \ , \lm_7=0 \ , \lm_8=1 \ . $$
\\
 And the normal lift is given by setting:
$$h_1=\gm_1 \ , \ h_2=\gm_2 \ , \\
h_3=1 \ , \ h_4=\lm_4 \ , \ h_5=\lm_5 \ , \ 
h_6=\lm_6 \ , \ h_7=\lm_7 \ , \ h_8=\lm_8 \ . $$
\medspace
We let the reader check that these lifts do satisfy the abnormal and normal
equations respectively.

\section{ Smoothness of Geodesics in the 2-Step Case}

In this section, we assume that $G$ is a nilpotent
lie group of 2-step. This  means that $[$\g$,[$\g,\g$]]=0$.
As before we consider a Carnot metric given
on a left-invariant distribution $D$  such 
that $D\oplus [$\g,\g$]=$\g. 

\begin{theorem}
\label{Tsmooth}
Under the above assumption, any  minimizer through $0$ is normal in
some subgroup of $G$ (for the induced Carnot metric) and hence any minimizer is smooth.
\end{theorem}
\remark
In the case where $G$ is a {\it free} nilpotent Lie group of 2-step,
one can prove that any minimizer is in fact normal,
see Gaveau \cite{Gav} and Brockett \cite{Br}.
\medspace
\Bproof
We need to show that any abnormal minimizer through $0$ is normal in some
subgroup of $G$. We will proceed by induction on the dimension
of $G$. The main step of the induction is given by the following lemma,
whose proof we postpone.

\begin{lemma}
\label{proper}
Any abnormal curve through $0$ (if it exists)
is tangent to a  left invariant proper sub-distribution 
$K\subset D$ and hence belongs to the proper Lie subgroup $H$ generated
by the algebra $K\oplus [K,K]$.
\end{lemma}

To start the induction, note that
any nilpotent Lie group of dimension 1 or 2 is in fact abelian,
and D must equal \g. Therefore there are no abnormal curves,
so every minimizer is normal 
(it is also easy to check that for dimension 3,
the only nilpotent Lie group is the Heisenberg group, which has no 
abnormal curves either).

Let $x(t)$ be an abnormal minimizer with $x(0)=0$. Let $H,K$ be as in 
Lemma \ref{proper} . The induced Carnot metric on $H$ is given by
the induced metric on $K\subset D$ (which as before we left-translate).
Since $x$ is a minimizer in $G$, it is a minimizer in $H$ (a horizontal
curve of smaller length than $x$ in $H$ would be horizontal and of smaller
length in $G$). Since $\dim H < \dim G$, $x(t)$ is normal, by induction.
If $x$ is a minimizer which does not pass through $0$, it is as smooth
as the minimizer $y(t)=x^{-1}(0)x(t)$ which does pass through $0$. \Eproof
\medspace
We now prove Lemma \ref{proper}.
\\
Because \g \ is 2-step, the equations for the abnormal curves 
(\ref{Eabnormal}) simplify to 
$$
\sum_{j=1}^r\sum_{k=r+1}^n\al_{ijk}\gm_j\lm_k = 0 \ \ \ \ {\rm for } \ \ 
i= 1,\ldots, r
$$
$$
\lm_k=0 \ \ \ \ {\rm for} \ \ k=1,\ldots,r \ \ \ \ {\rm and} \ \ 
\lm_k'= 0 \ \ \ \ {\rm for} \ \ k=r+1,\ldots,n 
$$
Indeed, in this case, $\al_{ijk}=0$ whenever $i$ or $j$ is greater
than $r$. 

Let $(x(t),\lm)$ be the cotangent lift of $x(t)$.
In particular  $x'=\gm$ satisfies the above equations, and $\lm\neq 0$.
\\
Let $M$ be the matrix whose entries are given by:
  
$$ m_{ij}=\sum_{k=r+1}^n \al_{ijk}\lm_k \ \ , \ i,j=1,\ldots,r $$

\begin{lemma}
\label{Lkernel}
If $\lm\neq 0$, the space $K={\rm Ker \ } M$ is a proper subspace of 
$\bbr^r$. Therefore \ $\{\sum \gm_j e_j \ | \ \gm \in K \}$
is a proper subspace of D.
\end{lemma}
Since for an abnormal minimizer $\lm \neq 0$, 
this claim obviously implies Lemma \ref{proper},
since  $x' \in {\rm Ker\ } M$.
\\
\begin{proof}
The fact that \g \ is 2-step (that is $D \oplus [D,D]=$\g)
implies the existence of coefficients
$\beta_{ijl}$ such that:
$$
e_l =\sum_{i,j\leq r}\beta_{ijl}[e_i,e_j], \ \ \ \ 
{\rm for}\ \ l= r+1,\ldots, n.
$$
{\bf Claim:} \ Viewing $A=\{\al_{ijk}\}$ and 
$B=\{\beta_{ijl}\}$ as $r^2\times (n-r)$ matrices 
(with $i,j$ ordered as a single index), the matrix $B^T$ is
a left inverse to $A$.
\ 
Indeed
$$
e_l =\sum_{i,j\leq r}\beta_{ijl}[e_i,e_j]= \sum_{i,j\leq r}\beta_{ijl}
\sum_{k>r}\al_{ijk}e_k=\sum_{k>r}
\left(\sum_{i,j\leq r}\al_{ijk}\beta_{ijl}\right)e_k,
$$
which implies
$$
\sum_{i,j\leq r}\al_{ijk}\beta_{ijl}= \delta_{kl}
\ \ \ ({\rm the \ Kronecker \ symbol}) \ , 
\ \ \ \ {\rm i.e.} \  A^TB= I_{n-r} \ .
$$
\medspace
The proof of Lemma \ref{Lkernel} is now easy. We want to prove that, unless
$\lm=0$, the matrix
$M$ is non zero. Suppose it were, i.e. $\sum_{k>r}\lm_k\al_{ijk}=0$, then
$$
\left(\sum_{k>r}\lm_k\al_{ijk}\right)\beta_{ijl}=0, \ \ \ \forall i,j,l
$$
\ 
$$
\Longrightarrow \ \ \ \  
0=\sum_{k>r} \ \sum_{i,j\leq r}\lm_k\al_{ijk}\beta_{ijl}=
\sum_{k>r}\lm_k\delta_{kl}=\lm_l
\ \ \ \forall l \ .
$$
\end{proof}

\ \bigspace
Christophe Gol\'e, Department of Mathematics, SUNY at Stony Brook,
Stony Brook, NY 11794. \\ gole@math.sunysb.edu
\medspace
Ron Karidi, Department of Mathematics, Stanford University,
Stanford, CA 94305. \\ karidi@math.stanford.edu

\end{document}